\documentclass[a4paper,12pt]{article}
\usepackage[utf8x]{inputenc}
\usepackage{amssymb, amsthm, amsmath, graphicx, fancyhdr, indentfirst, newlfont, latexsym, enumitem, wasysym,  textcomp}
\theoremstyle{remark} 
\newtheorem*{rem}{Remark}
\theoremstyle{plain} 
 \newtheorem{prop}{Proposition}
\theoremstyle{plain}
\newtheorem{teor}{Theorem}
\newtheorem*{teor*}{Theorem}
\theoremstyle{definition}

\newtheorem*{cor}{Corollary}
\newtheorem{defin}{Definition}

\newcommand{\hg}[1]{\nabla #1} 
 
\newcommand{\derpa}[3]{\dfrac{{\partial}^{#1}  #2}{ \partial {#3}^{#1}}}
\newcommand{\funzione}[5]{\begin{array}{cclcl} #1: & #2 & \longrightarrow & #3 \\ 
					           & #4 & \mapsto     & #5 \end{array}}

\newcommand{\dercov}[2]{\dfrac{D}{d #1}  #2}
\newcommand{\bb}{b }
\newcommand{\dbb}{\dot b}
\newcommand{\fifi}{\phi}
\newcommand{\dfifi}{\dot \phi}

\newcommand{\CVETT}[1]{\mathfrak{X} (#1)}

\newcommand{\HA}{H_{\alpha}^h}
\newcommand{\HB}{H_{\beta}^h}

\newcommand{\hp}{\mathbb{H}^2}

\newcommand{\N}{\mathbb{N}}
\newcommand{\roeucl}{\rho^{\R^2}}
\newcommand{\R}{\mathbb{R}}
\newcommand{\uah}{u_{\alpha}^h}

\newcommand{\CC}{C^{2,\delta}}

\newcommand{\Norm}{\eta}
\newcommand{\disco}[2]{B_{#1} \left( #2 \right)}
\newcommand{\sfera}[2]{S_{#1} \left( #2 \right)}

 \title{Constant mean curvature graphs on exterior domains of the hyperbolic plane.}
 \author{G. Citti and C. Senni}
\date{February, 2011}

\begin{document}

 \maketitle

 \begin{abstract}
We prove an existence result for non rotational constant mean curvature ends in $\hp \times \R$, where $\hp$ is the hyperbolic real plane. The value of the curvature is $h \, \in \, (0, \frac{1}{2})$. We use Schauder theory and a continuity method for solution of the prescribed mean curvature equation of exterior domains of $\hp$. We also prove a fine property of the asymptotic behavior of the rotational ends introduced by Sa Earp and Toubiana.
 \end{abstract}

\section*{Introduction}
The problem of ends of constant mean curvature, \textit{cmc} for short, in 3-manifolds is very classical and well understood in the Euclidean setting, at least in the finite total curvature case. In the 80's Schoen \cite{BELLARIC} proved that a finite total curvature minimal end is asymptotically a plane or a catenoid. This fact naturally arises the converse question: what curve can be a boundary of an end? A natural way to address the problem of existence of minimal ends is to look for solutions of the Dirichlet problem of constant zero mean curvature on an exterior domain. Hence one has to deal with the minimal surfaces equation on non convex domains, and it is well known (\cite{FI} and \cite{JS}) that non convexity can lead to non existence of the solution. This Dirichlet problem on non convex domains has been considered by Krust in \cite{KRUST}, Kuwert in \cite{KUWERT}, Tomi and Ye in \cite{TOMIYEYE}, Ripoll and Sauer in \cite{RIPOLLSAUER}. In the 90's, Kutev and Tomi solved the problem in \cite{KT1} and \cite{KT2} providing sufficient conditions on the boundary for the existence of minimal exterior graph with finite total  curvature.\\
The study of cmc ends in $\hp \, \times \, \R$ is at its first steps, but it turns out to be very interesting, since its properties seem to be quite different from the corresponding ones in $\R^3$.\\

In the Euclidean case it is possible to reduce the study of constant mean curvature, \textit{cmc} for short, surfaces to the two cases of zero and positive curvature, while in $\hp \, \times \, \R$ one has to distinguish at least three instances according to the value $h$ of the mean curvature: the minimal case, the case $h \, \in (0, \frac{1}{2}]$ and the case $h \, \in \, (\frac{1}{2}, +\infty)$. The role of the value $h=\frac{1}{2}$ has been outlined by Daniel on one side and Spruck for a different aspect. The first phenomenon occurring for $h = \frac{1}{2}$, discovered in \cite{BENOIT}, is the existence of a local isometry between minimal surfaces of the (Riemannian) Heisenberg group and surfaces of $\hp \times \R$ with constant mean curvature equal to $\frac{1}{2}$. The other one, described by Spruck in \cite{IGSP}, is that $h=\frac{1}{2}$ is the biggest value of the mean curvature such that the horosphere convexity of the boundary (see \cite{CIAO} and \cite{CM} for more details on this concept) is sufficient to have a solution of the constant mean curvature equation on a regular domain with any prescribed regular boundary value.\\
To study the problem of ends, one should consider the papers by Nelli and Rosenberg \cite{NEROSIN}, Sa Earp and Toubiana \cite{SATOU}, Spruck \cite{IGSP}, and very recently Elbert, Nelli and Sa Earp \cite{NESA},. In the first work it is proved that, for mean curvature $h$ between $0$ and $\frac{1}{2}$, there are no closed $h-$surfaces, hence it is natural to look for ends within these values of mean curvature. In this paper it is also proven that any $h-$graph on an exterior domain has non bounded height function. In the second paper, for each $h \, \in \, \left( 0, \frac{1}{2} \right]$, is introduced a one parameter family $\{ \HA \}_{\alpha}$  of rotational ends which have features similar to the ones of Euclidean catenoids. The problem of existence of non rotational ends has been considered by Nelli and Sa Earp in \cite{NESA}. \\ 
Our work deals with the $h \, \in \, (0, \frac{1}{2})$ case, and our main result is the existence of non rotational $h-$ends, Theorem \ref{teor: existence on exterior}. To prove our Theorem we require that the boundary of the exterior domain satisfies some geometric hypotheses. Our proof cannot be a simple adaptation of the analogous result in \cite{NESA} for the $h = \frac{1}{2}$ case because the asymptotic behavior of the $\HA$ functions in the two situations is quite different. This is why we establish a fine estimate (Theorem \ref{teor: dependance of deriv H in alfa}) clarifing the dependence on $\alpha$ of the asymptotic behavior of the $\HA$ family. \\
To obtain the existence result we use Schauder theory. First of all we establish a-priori estimates for a class of $h-$graphs on compact annuli, these estimates being done by means of the $\{ \HA \}_{\alpha}$ and of a perturbed distance function. After that, with a slightly modified method of continuity, we prove the existence of an $h-$graph on such a compact annulus. To obtain the vertical end, which is an end with unbounded height function, we consider a sequence of $h-$graphs on annuli diverging to the exterior domain and we prove that the estimates obtained in the compact cases are bounded on the sequence. As we have mentioned, the a-priori estimates are obtained by means of two results owing to the geometry of $\hp \, \times \, \R$. In particular to prove the estimates in a neighborhood of the inner boundary of the annulus we need to study the relation between the Laplacian of the distance from a regular compact set and the curvature of the level sets of this distance. In other words we need to understand the evolution of the curvature of a smooth Jordan curve along the flow of the distance function. In order to build barriers and uniform estimates on the exterior boundary, we need to establish a precise estimate of the asymptotic behavior of the $\{ \HA \}_{\alpha}$ family in dependence on the parameter $\alpha$:  we prove that, far enough from $0$, the function $\HA$ is monotonically decreasing on $\alpha$ in an open interval containing $(0, 2h]$. The proof of this fact is the most technical part of this work.\\
The contents are organized as follows:\\
In section one we recall a few facts of hyperbolic geometry.\\
In section two we study some geometry of $\hp \, \times \, \R$:  we prove the fine property of the asymptotic behavior, Theorem \ref{teor: dependance of deriv H in alfa}, and the result regarding the Laplacian of the distance. Here we describe the geometric hypotheses assumed in the existence Theorem \ref{teor: existence on exterior}. \\
In section three we present the a-priori estimates on compact annuli and the existence Theorems.

\section{Hyperbolic setting}
Here we recall only the properties of hyperbolic geometry we use in the paper, for an introduction to hyperbolic geometry we refer to \cite{ANDERSON}.\\
We consider the Poincar\'e's model of the hyperbolic plane, which means  
$$\hp = \{ z =(x,y) \, \in \R^2 : |z|_{\R^2} < 1\}$$ 
with the conformal metric 
\begin{align}
d \sigma^2(z) & = \left( \dfrac{2}{1- |z|_{\R^2}^2} \right)^2 \, \Big( dx^2 + dy^2 \Big) \nonumber \\
& =: \lambda^2(|z|_{\R^2}) \, \Big( dx^2 + dy^2 \Big) \label{eq: lambda}
\end{align}
 We denote by $\nabla^{\hp}$ the Levi Civita connection given by the metric and the curvature tensor by 
\begin{align}
 R(X,Y)Z = \nabla_{Y} \nabla_{X} \, Z - \nabla_{X} \nabla_{Y} \, Z + \nabla_{[X,Y]} \, Z \label{eq: curvature tensor}
\end{align}
We refer to $\partial \hp =  \{ |z|_{\R^2} =1\}$ as the \textit{asymptotic boundary} of $\hp$ because this is a set at infinite distance from any point of the hyperbolic plane.\\ In this model the geodesics are (suitable parametrizations of) arcs of Euclidean circles crossing orthogonally $\partial \hp$, or (suitable parametrizations of) Euclidean segments emanating from $0 \, \in \hp$. We also recall that by \textit{horocycle} are curves with constant geodesic curvature and equal to $1$. In our model horocycles are (suitable parametrizations of) euclidean circles tangent to $\partial \hp$.\\
$\hp$ is an homogeneous manifold with a three dimensional group of isometries and sectional curvature constant and equal to $-1$. The homogeneity allows us to identify any fixed point with $0 \, \in \, \hp$.\\
Let us recall that the isometry group of $\hp$ is generated by the three following elements
\begin{itemize}
 \item translations along geodesics
\item translation along horocycles
\item rotations about the point $0$
\end{itemize}
To describe the Riemannian product $\hp \, \times \, \R$ we use the coordinates given by the product. Assuming we use $t$ as a coordinate for $\R$, the metric we consider on $\hp \, \times \, \R$ is  $ ds^2 = d \sigma^2 + dt^2 $. In these coordinates if $S \, \subset \, \hp \times \R$ is a smooth surface that is a graph, i.e. $S = \{ (z, u(z)) : z \, \in \, \Omega\}$ for some $\Omega \, \subset \, \hp$ and $u \, \in \, C^2(\Omega)$, its mean curvature $H$ writes
\begin{align*}
  H(z) = - \dfrac{1}{2} \, \mathrm{div}_{\hp \, \times \, \R}(\eta(z,u(z))) = \dfrac{1}{2} \, \mathrm{div} \, \left( \dfrac{\nabla u(z)}{\sqrt{1 + {{|\nabla \, u(z)|}^2}}} \right)
\end{align*}
where $\eta$ is the upward unit normal vector to $S$ in $\hp \, \times \, \R$, $|\,.\,|$, $\nabla$ and $\mathrm{div}(.)$ are respectively the metric, connection and divergence of $\hp$. In this case the mean curvature can be considered as a second order differential operator that we will denote $Q$.  Hence  if $u \, \in \, C^2(\Omega)$ we denote
\begin{align*}
 Q(u) = \mathrm{div} \left( \dfrac{\nabla u}{\sqrt{1 + {{|\nabla \, u|}^2}}} \right)
\end{align*}
 It is well known that this operator, like the Euclidean one, is quasilinear and elliptic, and uniformly elliptic whenever $|\nabla u|$ is bounded, hence the theory of Schauder estimates and H\"older spaces can be used.

\subsection{The $\HA$ family}
From now on $h$ will be a real number belonging to the interval $(0, \frac{1}{2} )$. We recall the formulas and some relevant properties of the rotational cmc surfaces introduced by Sa Earp and Toubiana in \cite{SATOU}.   For $\alpha \, \in (0, +\infty)$ we denote
\begin{align}
\phi^h( \alpha) & =  \Big( \dfrac{-2 \, \alpha \, h + \sqrt{1-4h^2 + \alpha^2 }}{1 - 4h^2} \Big) \label{eq: phi nel raggio delle HA} \\
\rho^h(\alpha) & = \mathrm{arccosh} \left( \phi^h(\alpha) \right) \label{eq: rho^h_alpha}
\intertext{and $\forall \, \rho > \rho^h_{\alpha}$} 
u^h_{\alpha}(\rho) & = \dfrac{-\alpha + 2h \, \mathrm{cosh}(\rho)}{\sqrt{ \mathrm{sinh}(\rho)^2 - (-\alpha + 2h \, \mathrm{cosh}(\rho))^2}} \label{eq: u alfa}\\
\HA(\rho) & = \int_{\rho^h(\alpha)}^{\rho} \, \uah(w) \, d w \label{eq: acca alfa}
\end{align}
We remark that  $\rho^h(2h) = 0$ and then $H^h_{2h}$ is a simply connected entire graph which we denote $S^h$. If $\rho$ has the meaning of hyperbolic distance from $0 \, \in \, \hp$ these formulas define a family of rotational surfaces in $\hp \, \times \, \R$, where by rotational we mean invariant with respect to the rotation about the line $\{ z=0 \} \, \subset \, \hp \, \times \, \R$. Moreover each one of these surfaces is a graph defined in the complement of the disc $\disco{0}{\rho^h(\alpha)}$. The following proposition recalls some of the properties of the $\{ \HA \}_{\alpha}$ family. The proofs of the following statements can be found in \cite{SATOU}, \cite{NESASATOU}.

\begin{figure*}[h]
\begin{center}    
\includegraphics[width=12cm]{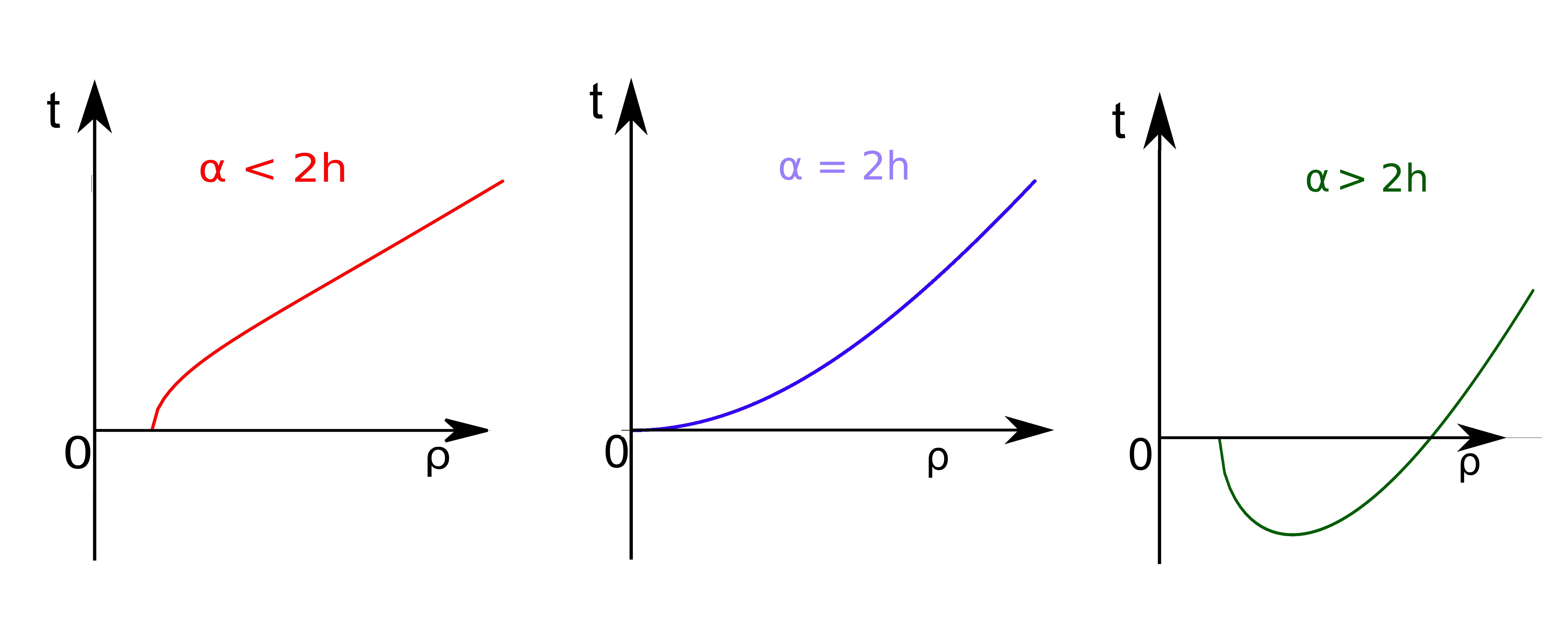}
\caption[legenda elenco figure]{Shape of the generating curves of $\{ \HA \}_{\alpha}$}\label{fig: HA_S_HB}
\end{center}
\end{figure*}

\begin{prop}\label{prop: prop of HA}$\\$
\begin{enumerate}
\item For all $ \alpha > 0 $ we have
$$
Q \Big( \HA \Big) \equiv 2 h
$$
\item If  $\alpha \neq 2h$, $\HA$ is zero valued and vertical on the hyperbolic circle $\sfera{0}{\rho^h(\alpha)}$
\item For $\alpha \leq 2h$ we have
\begin{align}
\HA(\rho) \geq 0 \quad \forall \, \rho \, \geq \, \rho^h(\alpha) \label{eq: HA non negativa}
\end{align}
For $\alpha \geq 2h$, $\HA$ is non-positive in a small annulus containing its boundary and positive out of this annulus. 
\end{enumerate}
\end{prop}
Figure \ref{fig: HA_S_HB} shows the dependence on the parameter $\alpha$ of the shape of the generating curves of $\{ \HA \}_{\alpha}$.

\section{Geometry of $\hp \, \times \, \R$}
We present two results we need to establish barriers for constant mean curvature graphs.\\
The first one deals with the flow of the distance function from a compact bounded by a regular Jordan curve. A standard technique to establish a-priori estimates for solutions of prescribed mean curvature problems (see for example \cite[Chapter 14]{GBT}) consists in bending graphs of suitable distance functions. Here we give an hyperbolic version of this technique and we choose as a distance function the distance from a compact set. Then we establish a sufficient condition on the curvature of the boundary so that the curvature of the level sets is a decreasing function of the distance.\\
The second result is about the asymptotic behavior of the family ${\{ \HA \}}_{\alpha}$. We prove that, for large $\rho$, the derivative of $\HA(\rho)$ with respect to $\alpha$ is negative. In other words we prove that the height $\HA(\rho)$ is monotonically decreasing in $\alpha$, at least in an open interval containing $(0, 2h]$.

\subsection{Evolution of the curvature along the flow of the distance function} 
We are interested to the evolution of the geodesic curvature of a regular closed curve  $\gamma \, \subset \, \hp$ along the flow of the distance function associated to the compact bounded by $\gamma$, thus we calculate the evolution of the geodesic curvature of the level sets associated with the flow of the field $\nabla d$.\\
Let be $q \in \hp$ and $t$ small enough, then by $\varphi_{q}(t)$ we denote the solution of the Cauchy problem
$$
 \left \{
\begin{array}{ccc}
\varphi_q(0) & = & q \\
\frac{d}{dt} \, \varphi_q(t) & =  & \nabla d (\varphi_q(t))
\end{array}	
\right.
$$

The following proposition gives explicit formulas for the evolution of the geodesic curvature along this flow.

\begin{prop}\label{prop: evoluz curv princ}$\\$
 Consider $\Omega \subset \hp$ a smooth domain bounded by $\gamma $ a Jordan curve and $d$ the distance function associated to $\Omega$, positive in $\Omega'$. Let be $q_0 \in \gamma$ and $k(q_0)$ the geodesic curvature of $\gamma$ in $q_0$. If we write $k_{q_0}(t)$ for the geodesic curvature of the curve $C_t = \{ p \, \in \Omega' : d(p) = t\}$ at the point $q(t)= \varphi_{q_0}(t)$ for $t$ small enough we have
\begin{align}
 |k_{q_0}|=1 & \Longrightarrow  k_{q_0}(t) \equiv  k_{q_0} \nonumber\\ 
|k_{q_0}|<1 & \Longrightarrow  k_{q_0}(t) = \tanh \left(t - \log \sqrt{\widetilde{k_{q_0}}} \right) \nonumber \\ 
|k_{q_0}|>1 & \Longrightarrow  k_{q_0}(t)= \coth \left(t - \log\sqrt{\widetilde{k_{q_0}}}\right)
\end{align}

where $\widetilde{k_{q_0}} = \left| \frac{1-k_{q_0}}{1+k_{q_0}}\right|$.
\proof
Consider $M^n$ an orientable complete manifold and $S \subset M^n$ a closed hypersurface. Denote by $\Omega$ the compact bounded by $S$ and, for $ a>0 $, denote $ V_a =  \{ p \, \in \, \Omega' : d(p) \leq a\}$. By compactness of $S$ we can use exponential coordinates on the tangent bundle of $M$ along $S$ and write $V_a = \{ (q,t)  \, : \, q \, \in \, S  \mbox{ and } 0 \leq t \leq a\} $. Moreover if
$$
S_t = \{ p \, \in \, \Omega' \, : \, d(p) = t\}
$$
for all $0 \leq t \leq a$ and $ q \in \, S$ we have
 $$
\hg{d}(q,t) \, \perp T_{(q,t)}S_t.
$$

 Then we can choose $\Norm_t = \hg{d}$ as a unit normal field along $S_t$, and  write $A^t$ for the shape operator associated with $\Norm_t$. If we consider $v_1(t), \dots v_{n-1}(t)$ an orthonormal frame of $T_{(q,t)} S_t$ diagonalizing $A^t$, the Radial Curvature Equation (see \cite[Theorem 3.6]{PETERSEN}, with signs adapted to our definition of the curvature tensor \eqref{eq: curvature tensor}) yields  
\begin{align}
 \left( - \nabla_{\hg{\eta_t}}A^t \right)v_i(t) + {\left( A^t \right)}^2 v_i(t) = R \left( v_i(t), \eta_t \right) \, \eta_t \qquad \forall \, i=1, \dots, n-1 \label{eq: rne on a frame} 
\end{align}
Since for all $ X \, \in \CVETT{S_t}$ we have
\begin{align*}
(\nabla_{\eta_t}A^t)X & = \nabla_{\eta_t}(A^tX) - A^t(\nabla_{\eta_t}X) \\
 & = \nabla_{\eta_t}(A^tX) \qquad \qquad \qquad \mbox{ because }  \nabla_{\eta_t}X \mbox{ is orthogonal to } S_t
\end{align*}
If $\frac{D}{dt}$ is the covariant derivative associated with the flow of $\eta_t$ we have
\begin{align*}
 \left( \nabla_{\eta_t }A^t \right) v_i(t) = \dercov{t}{\left( k_i(t) v_i(t)\right)}
\end{align*}
Thus, taking the scalar product with $v_i$ and using \eqref{eq: rne on a frame}, we get
\begin{align}
 -k_i'(t) + k_i(t)^2 = -\mathrm{sect}(v_i(t), \hg{d_t}) \qquad \forall \, i=1, \dots, n-1
\end{align}
In our case $M= \hp$ and $S=\gamma$, then we get $- \mathrm{sect}(v_1(t), \hg{d}) \equiv 1$ and for any $q_{0} \, \in \, \gamma$ the evolution equations of the curvature are:
$$
\left\{
\begin{array}{ccc}
  -k_{q_0}'(t) + k_{q_0}(t)^2 &=& 1 \\
k_{q_0}(0) &= & k_{q_0}
\end{array}
\right.
$$
Directly integrating we obtain the claim.

\endproof
\end{prop}

\begin{rem}
 \begin{itemize}$\\$
\item Following Cabezas-Rivas and Miquel \cite{CM}, we say that a curve is \textit{horosphere convex} when its geodesic curvature is greater than one.
\item The Proposition states in particular that if $|k_g(q)| \geq 1$ for all $q \in \gamma$, then the geodesic curvature is monotone decreasing along the flow. Thus if $\gamma$ is a horosphere convex curve in $\hp$, its geodesic curvature does not grow during the evolution along the distance flow but it stays greater than one.
\item If $\gamma$ is a circle, its evolution at any fixed time is a circle with bigger radius and hence smaller curvature.
 \end{itemize}
\end{rem}

\subsection{The asymptotic behavior of the $\HA$ family}
In this section we prove a crucial asymptotic property of the ${\{ \HA \}}_{\alpha}$ surfaces. We recall that Sa Earp and Toubiana in \cite{SATOU} introduced the explicit expression \eqref{eq: acca alfa} of $\HA(\rho)$. This expression together with Taylor approximation ensure that for large $\rho$
 \begin{align}\label{eq: sviluppo asintotico con alpha}
 H_{\alpha}^h(\rho) = k^h_{\alpha} + \dfrac{2h}{\sqrt{1-4h^2}} \, \rho +  O(e^{-\rho})
\end{align}
In this section we prove that the zero order term, $k^h_{\alpha}$, is strictly monotonically decreasing in $\alpha$, at least in an open interval containing $(0, 2h]$. 
Before proving the result we recall a simple and crucial property of the $\{ \HA \}_{\alpha}$ surfaces. By a straightforward calculation one can prove the following fact describing the dependence of the base circle of $\HA$ on the parameter.
\begin{prop}\label{prop: 1-1 alpha rho}$\\$
 \begin{itemize}
  \item $\rho^h(\alpha)$ is strictly monotonically decreasing on $(0, 2h]$ and $$\rho^h \Big( (0,2h] \Big) = \left[ 0, \dfrac{1}{\sqrt{1-4h^2}} \right)$$
\item $\rho^h(\alpha)$ is strictly monotonically increasing on $[2h, +\infty)$ and $$\rho^h \Big( [2h, +\infty) \Big) = [0, +\infty)$$
 \end{itemize}

\end{prop}

\begin{rem}$\\$
The behavior of the elements of ${\{ \HA \}}_{\alpha}$ near the base circle changes if $\alpha \, \in [0, 2h)$ or if $(2h, + \infty)$ (see figure \ref{fig: HA_S_HB}) and the same applies for the dependence of the radius of the base circle on the parameter. Hence it is useful to have two different notations for these two intervals of the parameter. We use the letter $\beta$ and we write $\HB$ when the parameter of the surface is greater than $2h$. Precisely when we write $\HB$, we tacitly assume $\beta > 2h$. 
\end{rem}
 We can now prove the Theorem describing the dependence on the parameter of the asymptotic expansion of $\HA$.

\begin{teor}\label{teor: dependance of deriv H in alfa}$\\$
 Let be $\rho$ large. Then  exists $ \overline \beta > 2h$ such that for all $ 0 < \widetilde \alpha \leq \overline \beta$ we have
\begin{align*}
 {\derpa{}{H^h_{\alpha}}{\, \alpha}(\rho)}_{| \, \alpha = \widetilde \alpha} < 0 
\end{align*}
\proof
This proof is made of two parts: first of all we prove that the derivative is negative for $\alpha < 2h$, then we prove that it blows up to $- \infty$ for $\alpha =  2h$.\\
 What we are going to do is a change of variable in $\uah$ (see  \eqref{eq: u alfa}) so that the derivation with respect to $\alpha$ does not interact with the singularity in $\rho^h(\alpha)$. Changing variable $s = \mathrm{cosh}(r)$, omitting some of the dependences on $\alpha$, we obtain:
\begin{align*}
 \HA(\rho) & = \int_{\fifi}^{\mathrm{cosh}(\rho)} \dfrac{-\alpha + 2h \, s}{\left( s-\bb \right)^{\frac{1}{2}} \left( s- \fifi\right)^{\frac{1}{2}}} \, \dfrac{ds}{\sqrt{s^2-1}}
\intertext{where}
\bb & = -\dfrac{2h \, \alpha \, + \, \sqrt{1-4h^2 + \alpha^2}}{(1-4h^2)}
\end{align*}
 is the negative zero of the denominator of $\uah$ and $\phi$ is defined in \eqref{eq: phi nel raggio delle HA}.\\
 To get rid of the dependence of the singularity on $\alpha$ we define $z = s-\fifi$ and $z(\rho, \alpha) = \mathrm{cosh}(\rho) - \fifi$ and we get
\begin{align}
 \HA(\rho)  & = \int_{0}^{z(\rho, \alpha)} \dfrac{-\alpha + 2h \, (z+\fifi)}{\sqrt{z} \, \sqrt{z + \fifi - \bb} \, \sqrt{{(z+\fifi)}^2-1}} \, dz \nonumber \\
\intertext{We denote }
    \widetilde u^h_{\alpha}(z)  & =: \dfrac{-\alpha + 2h \, (z+\fifi)}{\sqrt{z} \, \sqrt{z + \fifi - \bb} \, \sqrt{{(z+\fifi)}^2-1}} \label{eq: u_h tilde}
\end{align}
 We remark that $\fifi > 1$ and $-\bb >0$. Thus the only singularity of the integrand function in the interval  $[0, \mathrm{cosh}(\rho)-\fifi]$ is $0$.\\
Let's now compute the derivative
\begin{align}
\derpa{}{H^h_{\alpha}}{\alpha}(\rho) = \int_0^{z(\rho, \alpha)} \derpa{}{\widetilde u_{\alpha}^h}{\alpha}(z) \, dz \quad + \quad \widetilde u_{\alpha}^h(z(\rho, \alpha)) \, \derpa{}{z}{\alpha}(\rho, \alpha) \label{eq: derivata rispetto ad alfa come somma di due pezzi}
\end{align}
where, being $\rho$ large, the second term can be neglected because $\derpa{}{z}{\alpha}(\rho, \alpha)$ does not depend on $\rho$ and
\begin{align*}
 \widetilde u_{\alpha}^h(z(\rho, \alpha)) = \dfrac{1}{\mathrm{cosh}(\rho)} \, \dfrac{2h - \dfrac{\alpha}{\mathrm{cosh}(\rho)}}{\sqrt{1- \dfrac{\fifi}{\mathrm{cosh}(\rho)}} \, \sqrt{1- \dfrac{\bb}{\mathrm{cosh}(\rho)}} \, \sqrt{1- \dfrac{1}{\mathrm{cosh}(\rho)}}}
\end{align*}
is quickly decaying to zero as $\rho$ is growing to $\infty$. We now prove that the first term of \eqref{eq: derivata rispetto ad alfa come somma di due pezzi} is negative when $\rho$ is large.\\
To write the derivative of $\widetilde u_\alpha^h$ the following notation is useful (we use the \textit{dotted} notation for derivatives with respect to $\alpha$):
\begin{align*}
l & = z + \fifi \\
 \psi_1(z) & = \dfrac{-1 + 2h \, \dfifi }{{(z + \fifi- \bb)}^{\frac{1}{2}}{\left( {(z + \fifi)}^{2} - 1 \right) }^\frac{1}{2}} \\
 & =  \dfrac{-1 + 2h \, \dfifi }{{(l  - \bb)}^{\frac{1}{2}}{\left( {l}^{2} - 1 \right) }^\frac{1}{2}} \\
\\
\psi_2(z)  & = - \, (\dfifi - \dbb) \,  \dfrac{-\alpha + 2h \, (z + \fifi)}{2 \, {(z + \fifi - \bb)}^\frac{3}{2} \, {((z + \fifi)^2 -1)}^\frac{1}{2}}  \\
&=  - \, (\dfifi - \dbb) \,  \dfrac{-\alpha + 2h \, l}{2 \, {(l - \bb)}^\frac{3}{2} \, {(l\,^2 -1)}^\frac{1}{2}} \\
\\
\psi_3(z) & = - \dfifi \, \dfrac{(z + \fifi) \, (-\alpha + 2h \, (z + \fifi))}{ {(z + \fifi- \bb)}^\frac{1}{2} \, {((z + \fifi)^2 -1)}^\frac{3}{2}} \\
& =  - \dfifi \, \dfrac{l \, (-\alpha + 2h \, l )}{ {(l - \bb)}^\frac{1}{2} \, {(l ^2 -1)}^\frac{3}{2}}
\end{align*}
Thus we obtain:
\begin{align}\label{eq: derivata integranda H}
 \derpa{}{\widetilde u_{\alpha}^h}{\alpha}(z) = \dfrac{1}{\sqrt{z}} \, \Big( \psi_1(z) + \psi_2(z)+ \psi_3(z) \Big )
\end{align}
Many of the terms contained in the expressions of the $\psi_i$ have a sign which does not depend on $h, \, \alpha$ and $z$. Indeed for all $z \geq 0$ and $\alpha \neq 2h$ we have:
\begin{itemize}
 \item $z+ \fifi -1 \stackrel{h, \alpha}{>} 0$ being $\fifi = \mathrm{cosh}(\rho^h(\alpha))$ and $\rho^h(\alpha) = 0 \Leftrightarrow \alpha = 2h $
\item $z + \fifi - \bb \stackrel{h, \alpha}{>} 0$ being $\fifi - \bb = \dfrac{2 \, \sqrt{1-4h^2+\alpha^2}}{1-4h^2 }$
\item $ -1 + 2h \, \dfifi\stackrel{h, \alpha} < 0$ being $\dfifi= \dfrac{1}{1-4h^2} \, \Big( -2h + \dfrac{\alpha}{\sqrt{1-4h^2 + \alpha^2}}\Big)$ 
\item $\dbb < 0$ being $\dbb =  \dfrac{1}{1-4h^2} \, \Big( -2h - \dfrac{\alpha}{\sqrt{1-4h^2 + \alpha^2}}\Big)$
\item $\dot \fifi - \dot b > 0$ being $\fifi - \bb = \dfrac{2 \, \sqrt{1-4h^2+\alpha^2}}{1-4h^2 }$
\end{itemize}

If moreover $\alpha <  2h$ we also have 
$-\alpha + 2h \, (z + \fifi) \stackrel{h, \alpha}{>} 0 $  which implies
$$
\psi_2(z) \stackrel{h, \alpha}{<} 0
$$

To prove that the first term of equation \eqref{eq: derivata rispetto ad alfa come somma di due pezzi} is negative it is thus enough to prove that 
$$ - \dfrac{\psi_1(z)}{\psi_3(z)} \geq 1$$ 
To simplify calculations we denote
\begin{align*}
 c_1 = - (-1 + 2h \, \dfifi) > 0 \quad c_{31}= - 2h \, \dfifi> 0\quad c_{32}- \alpha \dfifi>0 
\end{align*}
and thus we get
\begin{align*}
-\dfrac{\psi_1(z)}{\psi_3(z)} = c_1 \, \dfrac{l^2 -1 }{c_{31}l^2 - c_{32}l} \geq 1 \Longleftrightarrow (c_1 - c_{31}) l^2 - c_{32} l - c_1 \geq 0
\end{align*}
where 
\begin{align*} 
c_1 - c_{31} = 1 
\end{align*}
The solutions of the equation
$$
l^2  + c_{32} l - c_1 =0
$$
are the two distinct real numbers
$$
l_{\pm} = \dfrac{- c_{32} \pm \sqrt{{c_{32}}^2 + 4 \, c_1 }}{2}
$$
A straightforward computation shows that $\mathrm{max} \lbrace l_{-}, \, l_{+} \rbrace \leq \phi$ and hence that, when $\rho$ is big enough and $\alpha \, \in (0,2h)$, $ \HA(\rho)$ is strictly decreasing when $\alpha$ increases.\\
If we evaluate equation \eqref{eq: derivata integranda H} in $\alpha = 2h$ we obtain a singularity in $0$ which is non integrable. Indeed being $ {\psi_3}_{|\alpha = 2h} \equiv 0$ we get 
\begin{align*}
\psi_1(z) & = - \dfrac{1}{\sqrt{z}} \, \dfrac{1}{\Big( z + \dfrac{2}{1-4h^2} \Big)^{\frac{1}{2}} \, \Big( z + 2 \Big)^{\frac{1}{2}}} \\
\psi_2(z) & = - \dfrac{1}{\sqrt{z}} \, \dfrac{4h^2 \, z}{(1-4h^2) \, \Big( z + \dfrac{2}{1-4h^2} \Big)^{\frac{3}{2}} \, \Big( z + 2 \Big)^{\frac{1}{2}}} \\
\intertext{and hence}
 \psi_1(z) + \psi_2(z) & =  - \dfrac{1}{\sqrt{z}} \left( \dfrac{\sqrt{z+2}}{(1-4h^2) \, \left(z + \dfrac{2}{1-4h^2} \right)^{\frac{3}{2}}} \right)
\end{align*}
So we can write 
\begin{align*}
 \derpa{}{\widetilde u_{\alpha}}{\alpha}(z)_{|\alpha=2h} & = \dfrac{\psi_1(z) + \psi_2(z) + \psi_3(z)}{\sqrt{z}} \\
& = - \dfrac{1}{z} \; \left( \dfrac{\sqrt{z+2}}{(1-4h^2) \, \left(z + \dfrac{2}{1-4h^2} \right)^{\frac{3}{2}}} \right)
\end{align*}
\endproof
\end{teor}
\begin{rem}$\\$
 \begin{itemize}
  \item This proof cannot be used in the case $\alpha > 2h$ because here $\widetilde u_{\alpha}^h(z)$ is positive near $z=0$. Indeed on can easily check that $\psi_1(0) = -\psi_3(0)$ and $\psi_2(0)>0$.
\item Roughly speaking we have proven that at infinity $\{\HA \}_{\alpha}$ is a family of cones with same angle, given by the mean curvature, and different vertex, given by the parameter and the curvature.
 \end{itemize}
\end{rem}
This Theorem allows us to associate two concepts of distance between rotational surfaces.
\begin{defin}{Asymptotic vertical distance}\label{def: asymptotic distance}$\\$
Consider $0 < \alpha_1 < \alpha_2 < \overline{\beta}$.
 We define the \textit{asymptotic vertical distance} between $H^h_{\alpha_1}$ and $H^h_{\alpha_2}$ to be 
$$
\lim_{\rho \to + \infty}  H^h_{\alpha_1}(\rho) - H^h_{\alpha_2}(\rho)
$$
\end{defin}

\begin{figure*}[h]
\begin{center}    
\includegraphics[width=12cm]{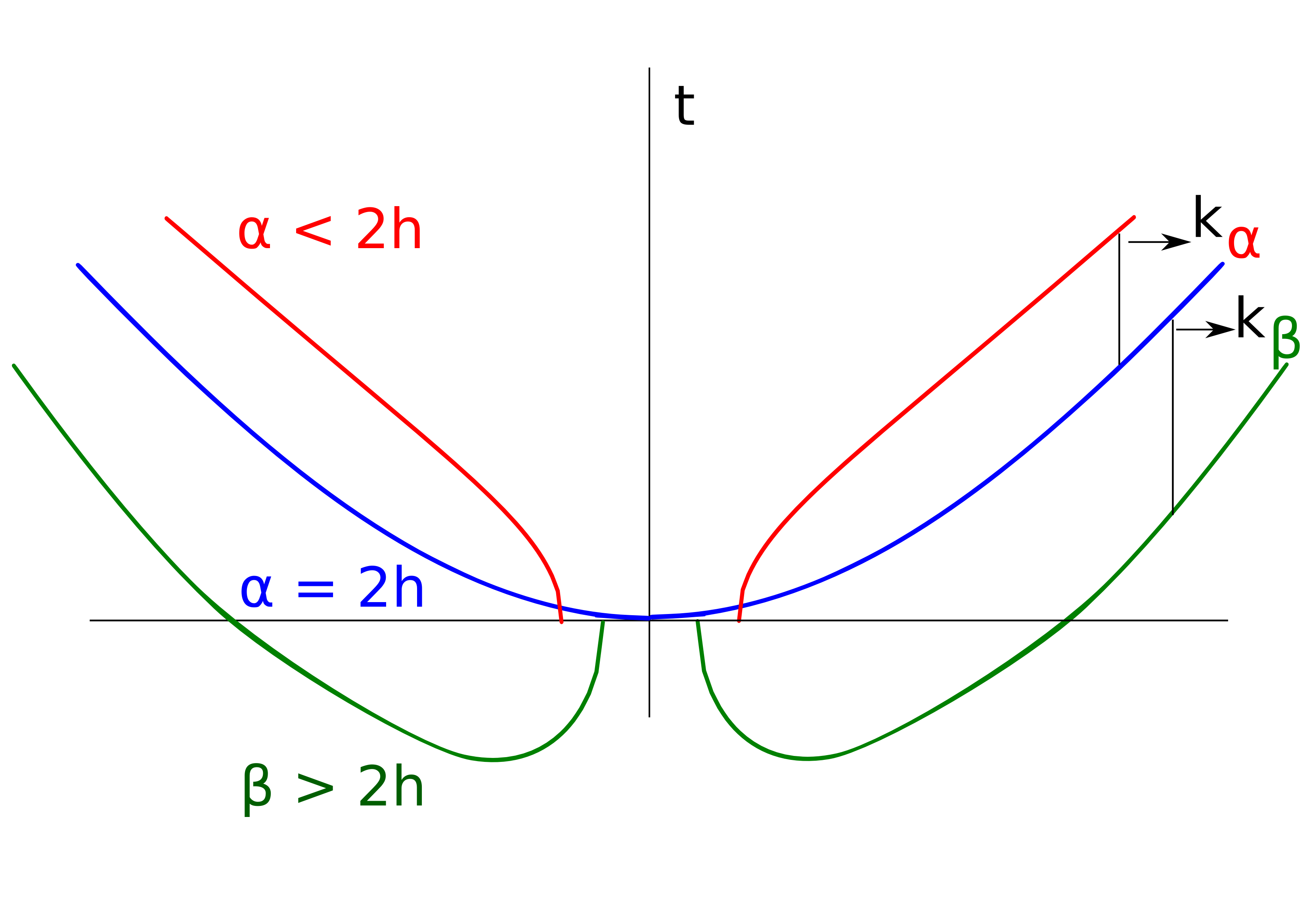}
\caption{Asymptotic behavior of $\{ \HA \}_{\alpha}$}\label{fig: sovrapposte}
\end{center}
\end{figure*}

Now we can explicit the consequences of the last Theorem in terms of the relative positions of the elements of the family ${\{ \HA \}}_{\alpha}$ according to the value of the parameter. 
\begin{cor}\label{prop: comportamento asintotico}$\\$
 Let   $h \, \in \, (0, \frac{1}{2})$. 
\begin{itemize}
\item For all $ \alpha \, \in \, [0, 2h) \,$
 $S^h \cap \HA$ is a circle. Inside this circle we have $\HA < S^h$, outside we have the opposite inequality.
\item Exists $ \overline \beta \, \in \, (2h, +\infty) \, $ such that  for all $ \beta \, \in \, (2h, \overline \beta)$ we have
\begin{align*}
 S^h & \cap \HB  = \emptyset \\
\HB & < S^h  \\
\HB & < \HA  \qquad \forall \, \alpha \in (0, 2h]
\end{align*}

\item The vertical asymptotic distance between $S^h$ and $H^h_{\alpha}$ is a positive real number if $\alpha > 2h$, negative id $\alpha<2h$. 

\end{itemize}
 \end{cor}
Figure \ref{fig: sovrapposte} shows the positions of elements of $\{ \HA \}_{\alpha}$ for different values of the parameter.

We end this section defining the concept of asymptotic horizontal distance for element of the ${\{ \HA \}}_{\alpha}$ family. This is the limit for $t \to +\infty$ of the distance between the circles obtained intersecting two distinct elements of the family with the plane of height $t$.

\begin{defin}{Asymptotic horizontal distance}$\\$
Consider $0 \leq \alpha \leq 2h$ and $2h < \beta \leq \overline \beta$. Since $\HA$ and $\HB$ are asymptotically invertible being their asymptotic behavior  affine, their inverse $\Big( \HB \Big)^{-1}$ and $\Big( \HA \Big)^{-1}$ are also affine. We define the \textit{asymptotic horizontal distance} between $\HA$ and $\HB$ to be
$$
d_{\infty}(\HA, \HB) = \lim_{t \to +\infty } \left| \Big( \HA \Big)^{-1}(t) - \Big( \HB \Big)^{-1}(t) \right|
$$
\end{defin}
Theorem \ref{teor: dependance of deriv H in alfa} ensures that both the asymptotic distances are non zero for $\alpha \neq \beta$.
\begin{prop}$\\$
Let be $0< \alpha, \beta < \overline \beta$. Then
\begin{align*}
 d_{\infty}(\HA, \HB)  = \left| \dfrac{k^h_{\beta} - k^h_{\alpha}}{c_h} \right| \label{eq: d_infty}
\end{align*}
where $c_h = \frac{2h}{\sqrt{1- 4h^2}}$
\proof

Denote by
$$ 
\rho_{\alpha,t} = \Big( \HA \Big)^{-1}(t)
$$ 
the radius of the circle in $\hp$ whose image via $\HA$ is contained in the plane $\hp \, \times \, \{ t\}$. \\
The asymptotic expansions \eqref{eq: sviluppo asintotico con alpha} of $\HA$ and $\HB$ yields
\begin{align}
\rho_{\alpha,t} - \rho_{\beta,t} & = \dfrac{k_{\beta}^h - k_{\alpha}^h}{c_h} + O(e^{-\rho})
\end{align}
\endproof
\end{prop}

\subsection{r-admissibility}
\begin{figure*}[h]
\begin{center}    
\includegraphics[width=12cm]{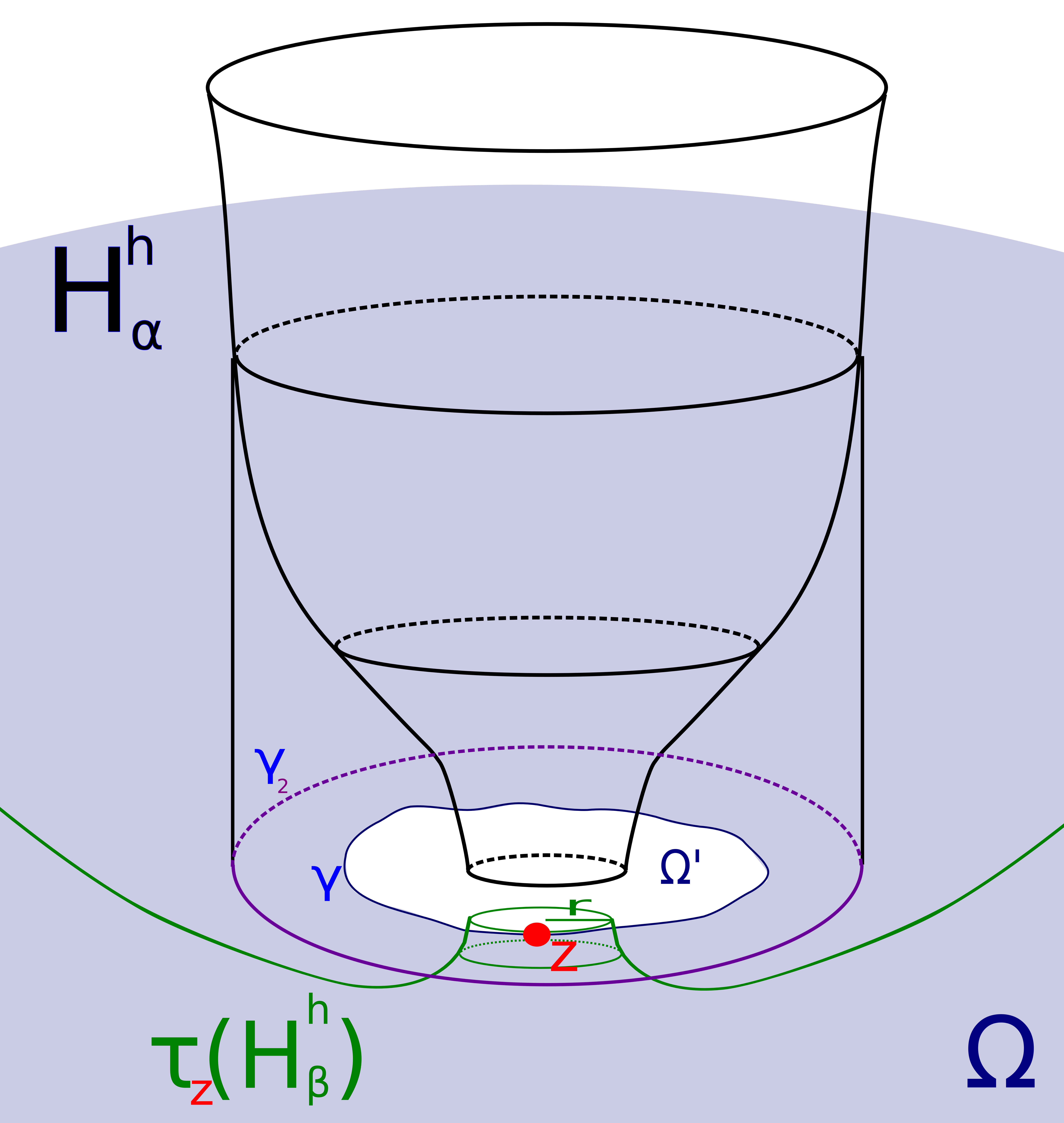}
\caption[legenda elenco figure]{$\HB$ as a barrier}\label{fig: sfera barriera}
\end{center}
\end{figure*}

Consider $\Omega \subset \hp$ a smooth annulus and $u : \Omega \to \R$ a function of class $C^2$ whose graph has constant mean curvature. As we have mentioned, we will use the $\{ \HA \}_{\alpha}$ surfaces to give  a-priori $C^1$ estimates of $u$. Our estimates hold for annuli satisfying geometric hypotheses, we now explain the one we called the $r-$admissibility hypothesis. Assume that the inner boundary of $\Omega$, $\gamma$, is a Jordan curve and the outer boundary, $\gamma_2$, is a circle centered in $0$. Let be $\alpha \, \in \, (0, 2h)$ such that $\sfera{0}{\rho^h(\alpha)} \subset \Omega$ and consider $u$ an $h-$graph on $\Omega$ satisfying $u_{|\gamma} = 0$ and $u_{| \gamma_2} = {\HA}_{| \gamma_2}$. Now take  $z \, \in \, \gamma$. We will consider a precise $\beta \, \in \, [2h, \overline{\beta}]$ depending on $\gamma$ and an horizontal translation $\tau_z$  making $\HB$ tangent to $\gamma$ in $z$. This $\tau_z(\HB)$ will be a lower barrier for $u$ in $z$ provided it is below $u$.
Figure \ref{fig: sfera barriera} shows how we will use $\HB$ as a barrier.\\

\begin{figure*}[h]
\begin{center}    
\includegraphics[width=12cm]{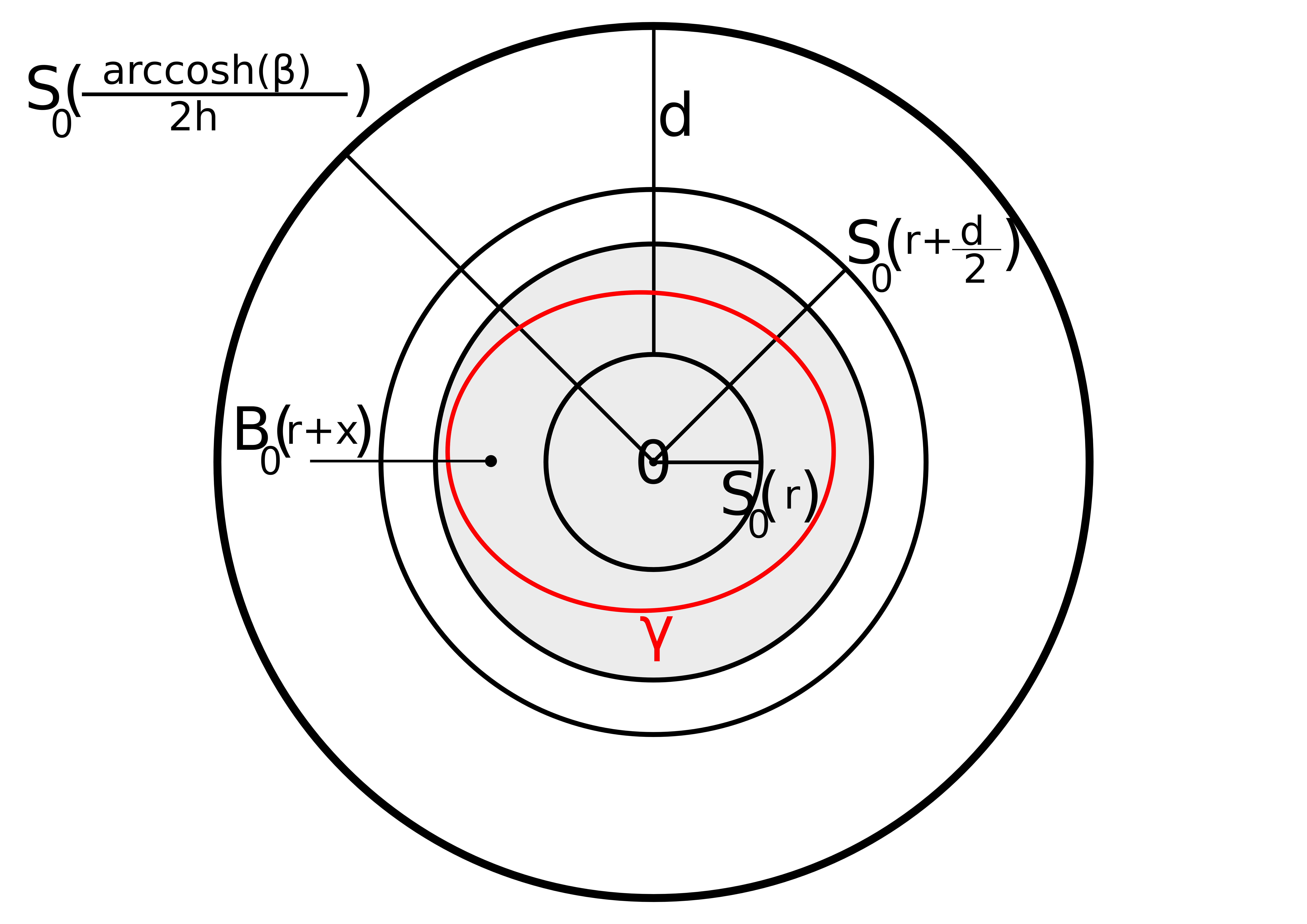}
\caption[legenda elenco figure]{$r-$admissible curve}\label{fig: curve e circonferenze traslazione}
\end{center}
\end{figure*}

 The aim of the  $r-$admissibility condition is to guarantee that $\tau_z(\HB)$ is actually below $u$. Moreover we have to guarantee that $\tau(\HB)_{|\gamma} \leq 0$.\\
 We now make these remarks precise. Consider a compact domain $\Omega' \, \subset \, \hp$ bounded by the smooth Jordan curve $\gamma$. Assume the curve satisfies an interior (with respect to $\Omega'$) sphere condition of radius $r >0$. We associate to this curve two elements of the $\{ \HA \}_{\alpha}$ family. First of all we consider a $\HA$ whose base circle is contained in the interior of $\Omega'$. More precisely we consider $\alpha \, \in \, (0, 2h)$ such that
\begin{align}
  H^h_{ \alpha} \cap \lbrace t=0 \rbrace \subset \, \mathrm{int}(\Omega') \label{eq: alpha che scelgo per l'esistenza} 
\end{align}
The other surface is given by the interior $r-$sphere condition as well: we take the $\beta \, \in \, (2h, +\infty)$ such that $\rho^h(\beta) = r$ if $\beta \leq \overline \beta$, where $\overline \beta$ is given in Theorem \ref{teor: dependance of deriv H in alfa}. If the value of $\beta$ such that $\rho^h(\beta) = r$ is bigger than $\overline \beta$, we reduce $r$ to $\widetilde r$ in such a way that $\beta(\widetilde r) \leq \overline \beta$. This is possible by proposition \ref{prop: 1-1 alpha rho}. Moreover, being $\widetilde r < r $, $\gamma$ also  satisfies an interior (with respect to the compact bounded by $\gamma$) $\widetilde{ r}-$sphere condition.\\
Now we define $d_{\beta}$ to be the hyperbolic distance between the base circle of $\HB$ and the circle where it has its minima (recall the definition of $\HB$ \eqref{eq: acca alfa}), i.e. the hyperbolic circle of radius $\frac{\mathrm{arccosh}(\beta)}{2 h}$. Thus we have 
$$
d_{\beta} = \dfrac{\mathrm{arccosh}(\beta)}{2 h} - \rho^h(\beta)
$$
By Theorem \ref{teor: dependance of deriv H in alfa} we can define the quantity
\begin{align}
 \xi = \min \Big \lbrace \dfrac{d_{\beta}}{2}, \; d_{\infty}(\alpha, \beta) \Big \rbrace > 0 \label{def: x}
\end{align}
where $d_{\infty}$ is defined in \eqref{eq: d_infty}. Remark that since $\alpha$ and $\beta$ depend only on $\gamma$, the same holds for $\xi$. We can now give the definition $r-$admissibility.
\begin{defin}{$r-$admissibility}\label{def: curva r-ammissibile}$\\$
 Let $\gamma \subset \hp$ be a smooth Jordan curve. Suppose $\gamma $ satisfies an interior sphere condition of radius $r$. We say that $\gamma$ is $r-$admissible if it is contained in the circular annulus 
\begin{align}
A_{\gamma} = \{ z \in \hp \, : \, \rho^h(\beta) \leq |z| \leq \rho^h(\beta) + \xi \} \label{eq: A_gamma}
\end{align}

\end{defin}
Figure \ref{fig: curve e circonferenze traslazione} describes $A_{\gamma}$.
 We now check that a curve verifying this definition allows the $\HB$ to be used as a $C^1$ barrier from below on $\gamma_1$.

\begin{prop}\label{prop: r-admissible implica barriera}$\\$
 Let $\gamma$ be an $r-$admissible curve and take $z \, \in \gamma$. Consider $\tau_z(\HB)$ an horizontal translation of $\HB$ tangent to $\gamma$ in $z$.\\
 Thus 
\begin{enumerate}
 \item  $\tau_z(\HB) < \HA$
\item $\tau_z(\HB)_{|{\gamma}} \leq 0$
\end{enumerate}
 
\proof$\\$
Recall that $\beta$ is chosen so that $\beta \in (2h, \bar \beta]$ and $\rho^h(\beta) = r$.\\ 
The first fact follows directly from the definition of $r-$admissibility, indeed the distance that has to be covered to make the circle $\sfera{0}{r}$ tangent to $\sfera{0}{r + \xi}$, precisely $\xi$, is not bigger than $d_{\infty}(\alpha, \beta)$ (recall \eqref{eq: d_infty} and \eqref{def: x}). Being $\gamma \subset A_{\gamma}$ by hypothesis of $r-$admissibility, the distance that needs to be covered to make $\sfera{0}{r}$ tangent to $\gamma$ is strictly smaller than the distance between the two circles, namely $\xi$. \\
To prove the second statement we move $\gamma$ instead of $\HB$. If we prove that any translation making $\sfera{0}{r + \frac{d}{2}}$ tangent to $\sfera{0}{r}$ is contained in the disc $\disco{0}{\frac{\mathrm{arccosh}(\beta)}{2h}}$ we are done. Indeed in this disc $\HB$ is negative and to make $\sfera{0}{\frac{d}{2}+ r}$ tangent to $\sfera{0}{r}$ we need to cover a distance greater than the distance we have to cover to make $\gamma$ tangent to $\sfera{0}{r}$. But $\sfera{0}{\frac{d}{2} + r} \, \subset \, \disco{0}{\mathrm{arccosh} \left( \frac{\beta}{2h} \right)}$.
\endproof
\end{prop}
Let's discuss some examples of $r-$admissible curves.
\begin{itemize}
 \item A circle of radius $r$ is $r-$admissible provided $r$ is small enough to have $\beta(r) \in (2h, \bar \beta]$ 
\item A more interesting example is given by a small $C^2$ perturbation of a small circle. Consider $\alpha \, \in [0, 2h)$ and $\beta \, \in \, (2h, \overline \beta]$ and denote 
\begin{align}
A(\alpha, \beta) = \{ z \, \in \, \hp : \, \rho^h(\beta) \leq |z| \leq \rho^h(\beta) + \xi \}
\end{align}
where $\xi$ is defined in \eqref{def: x}. Then, every circle contained in $A(\alpha, \beta)$ is $r-$admissible for $r = \rho^h(\beta)$. Moreover if we consider a perturbation of such a circle with small $C^0$ and $C^2$ norms, we still have an $r-$admissible curve. Indeed a small $C^0$ norm ensures that the deformed curve is still in $A(\alpha, \beta)$ and a small $C^2$ norm ensures we have haven't changed too much the geodesic curvature. Hence for the deformed curve we can still use the same interior sphere we used for the initial circle.\\
\item Another example of $r-$admissible curve is given by the constant volume mean curvature flow considered by Cabezas-Rivas and Miquel in \cite{CM}.  Indeed we have the following Corollary of the Theorem describing the evolution via constant volume mean curvature flow of a horosphere convex domain in $\hp$ (see \cite[Theorem 1.2]{CM}).
\begin{cor}\label{prop: perche il flusso cvmcf conserva la r-ammissibilita}$\\$
 Let $\gamma \, \subset \, \hp$ a smooth Jordan curve and consider $X$ its constant volume mean curvature flow and call $C$ its limit circle. Then for all $ \varepsilon >0$ exists $t(\varepsilon) \, \in \, \R^+$ such that for all $ t \geq t(\varepsilon)$ we have
\begin{align*}
 d(X(t), C) & \leq \varepsilon \\ 
|k_g(X(t)) - k_g(C) | & \leq \varepsilon
\end{align*}
\end{cor}
Thus the constant volume mean curvature flow is an useful tool to deform a curve not too different from a circle into a circle preserving $r-$admissibility.\\
 
We remark that all the examples just described can be deformed into a circle preserving horosphere convexity and $r-$admissibility. 
\end{itemize}

\section{Existence results}

In this section we prove the main result of this paper which is the existence of non rotational ends with constant mean curvature $h \, \in \, (0, \frac{1}{2})$. First of all we prove a-priori estimates for $h-$graphs on a class of compact annuli, then we prove existence  in the compact case and finally we prove existence of vertical ends as limit of a sequence of compact graphs. In this section $\Omega$ will be a compact smooth annulus with boundary two curves $\gamma_1$ and $\gamma_2$, where $\gamma_1$ is contained in the compact bounded by $\gamma_2$. All circles are thought as centered in $0 \, \in \, \hp$ and we also consider, without loosing any generality, annuli bounding a compact domain containing $0$. Before proving the results, we explicit that in all the proofs by $\HA$ we mean a small negative vertical translation of $\HA$ and by $\HB$ we  mean a small positive vertical translation of $\HB$. Precisely, making an abuse of language, we denote by $\HA$ the surface $\HA - \varepsilon$ with $0 < \varepsilon << 1$ and we denote by $\HB$ the surface $\HB + \varepsilon$ with $0 < \varepsilon << 1$. This is because we use these surfaces as barriers for the gradient on the boundaries and hence we need the normal derivative to be finite. We also explicit that we compute normal derivatives according to the inner normal, hence lower barriers will bound normal derivatives from below while upper barriers will bound normal derivatives from above.\\

\begin{teor}\label{teor: estimates on annuli}$\\$
Let be $h \, \in \, (0, \frac{1}{2})$ and let $\Omega \, \subset \, \hp$ be a compact smooth annulus. Assume $\gamma_1$ is $r-$ admissible, horosphere convex and $\gamma_2$ is a circle with a large radius. Consider $\alpha = \alpha(r)$ as in the definition of $r-$admissibility. If for some $0 < \delta < 1 $, $u \, \in C^{2,\delta}(\Omega)$ is a solution of the Dirichlet problem
\begin{equation} \label{eq: Dirichlet caso compatto}
\left \{
\begin{array}{rccr}
 \dfrac{1}{2} \, \mathrm{div} \left( \dfrac{\nabla u}{\sqrt{1 + {|\nabla u|}^2}} \right) & = & h & \mbox{ in }  \Omega \\
u & = & 0 & \mbox{ in }  \gamma_1 \tag{D} \\
u & = & H^h_{\alpha(r)} & \mbox{ in } \gamma_2
\end{array}
\right.
\end{equation}
then exists $ C=C(h, \gamma) > 0$ such that
\begin{equation}
 || u ||_{C^{2,\delta}(\Omega)} \leq C
\end{equation}
\proof$\\$
\begin{figure*}[ht]
\begin{center}    
\includegraphics[width=11cm]{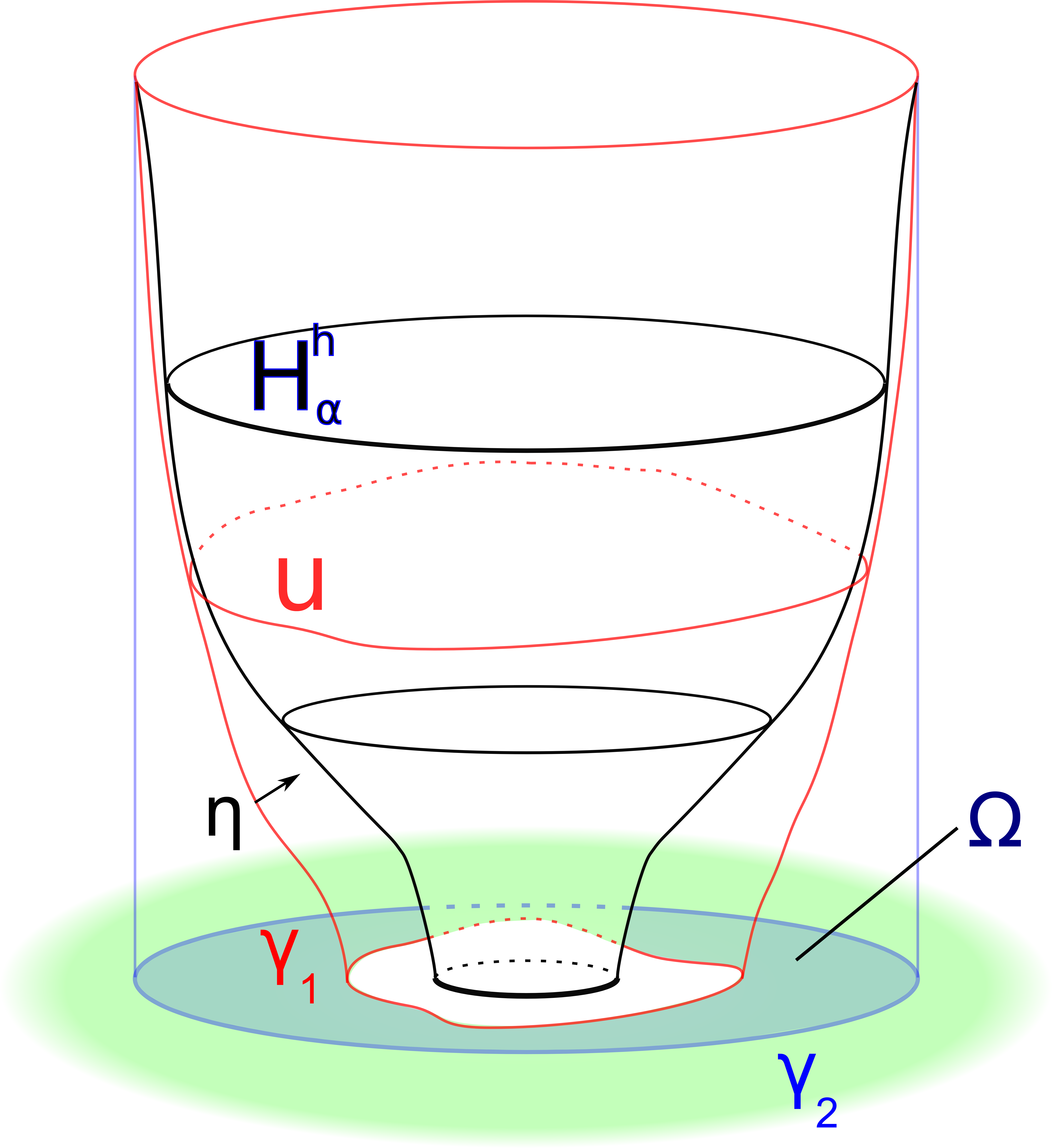}
\caption{A priori estimates}\label{fig: a priori estimates}
\end{center}
\end{figure*} 
Figure \ref{fig: a priori estimates} sketches the picture of the $h-$graph for which we are going to give the a-priori estimates.\\
By standard Schauder theory $||u ||_{C^1(\Omega)}$ estimates imply $|| u ||_{C^{2, \delta}(\Omega)}$ estimates, hence we are going to establish estimates of $u$, boundary gradient estimates of $u$ and interior gradient estimates of $u$.\\
Let be $\beta \, \in \, (2h, \overline{\beta} )$ given by the hypothesis of $r-$admissibility of $\gamma_1$.

\begin{itemize}
 \item \textit{$C^0$ estimates}\\ 
We use the maximum principle to prove
$$
\HB \leq u \leq \HA \mbox{ on } \Omega
$$
Let's prove that $u$ is below $\HA$:
\begin{align*}
 {u}_{|\gamma_1} = \; & 0 <  {\HA}_{|\gamma_1} & \mbox{  $\HA$ being positive on $\gamma_1$ by \eqref{eq: HA non negativa}} \\
{u}_{|\gamma_2}  = & \;  {\HA}_{|\gamma_2}  &  u \mbox{ being a solution of \eqref{eq: Dirichlet caso compatto} } \\ 
Q(u) = & 2h = Q(\HA) &
\end{align*}
and hence $u \leq \HA$ all over $\Omega$.\\
Let's prove that $u$ is above $\HB$:
\begin{align*}
 {u}_{|\gamma_1} = \; & 0  > {\HB}_{|\gamma_1} & \mbox{ by $r-$admissibility hypothesis} \\
{u}_{|\gamma_2} = & \; {\HA}_{|\gamma_2}  > {\HB}_{|\gamma_2} & \mbox{by Theorem \ref{teor: dependance of deriv H in alfa}} \\
Q(u) = & 2h = Q(\HB) & 
\end{align*}
and hence $u \geq \HB$ all over $\Omega$.\\
We have proved that the graph of $u$ is in the compact region contained between the graph of $\HA$ and the graph of $\HB$.

\item \textit{$C^1$ estimates on $\partial \, \Omega$ }\\
We will prove four estimates, using either $\HA$ and $\HB$ as barriers or a distance function.

\subitem \textit{$C^1$ estimates on $\gamma_1$ }\\ 
We build an upper barrier bending a distance function and we use the $\HB$ as a lower barrier.\\
Consider $z \, \in \, \Omega$ and define $d(z) = d_{\hp}(z, \gamma_1)$. 

Recall that for $\varepsilon > 0$ by $V_{\varepsilon}$ we mean the part of the tubular neighborhood of $\gamma_1$ of thickness $\varepsilon$ contained in $\Omega$. We are going to show that exist $\varepsilon > 0$ and $A>0$ such that the function $w^+ := \psi_{\varepsilon, A} \circ d$ defined by
\begin{align}
\funzione{\psi_{\varepsilon, A}}{[0, \varepsilon]}{\R}{x}{A \, \left( e^{\varepsilon} - x - e^{\varepsilon - x}\right)}
\end{align}
is an upper barrier for the normal derivative of $u$ on $\gamma_1$. First of all denote
\begin{align}
c= c(\gamma_1) =  \max_{y_0 \, \in \, \gamma_1} k_g(y_0) \label{eq: max geod curvature along gamma_1}
\end{align}

\begin{itemize}
 \item We are going to prove that there is an $\varepsilon >0$ such that
$$
Q(w^+) < 0 \mbox{ in } V_{\epsilon}
$$
where $Q$ is the mean curvature operator. For any $\varepsilon > 0$ by a simple calculation using general properties of distance functions one can see that for all $ q \, \in \, V_{\varepsilon}$
$$
\dfrac{1}{2} \, Q(w^+)(q) = \dfrac{1}{\sqrt{1+{\psi'}_{\varepsilon, A}^2}} \Big( {\psi'}_{\varepsilon, A} \, k_g(q) + \dfrac{{\psi}_{\varepsilon, A}''}{1+{\psi'}_{\varepsilon, A}^2}\Big)
$$
where $k_g(q)$ is the geodesic curvature of the curve $\{ p \in \Omega : d(p)=d(q) \}$ calculated with respect to the unit normal vector $\Norm =-\hg{d}$ in $q$. Using the horosphere convexity of $\gamma_1$ and Proposition \ref{prop: evoluz curv princ} we write
$$
Q(w^+)(q) \leq \dfrac{1}{\sqrt{1+{\psi'}_{\varepsilon, A}^2}} \left( \psi_{\varepsilon, A}'\, k_g(z) + \dfrac{\psi_{\varepsilon, A}''}{1+{\psi'}_{\varepsilon, A}^2}\right)
$$
where $q =\varphi_{z}(t)$ for a unique $z \, \in \, \gamma_1$ and $0\leq t \leq \varepsilon$ (recall that $\varphi_{z}(t)$ is the flow of $d$).  Thus we get
\begin{align}
 c \, \psi_{\varepsilon, A}' \, \Big( 1 + \psi_{\varepsilon, A}'^2 \Big) + \psi_{\varepsilon, A}'' \leq 0  \quad \Longrightarrow \quad  Q(w^+)(q) \leq 0
\end{align}
By a straightforward calculation one can verify that
\begin{align}
 \Big( c \, \psi_{\varepsilon, A}' \, \left( 1 + \psi_{\varepsilon, A}'^2 \right) + \psi_{\varepsilon, A}'' \Big)_{|\{ \varepsilon = 0, x= 0\}} = -A <0
\end{align}
hence there is a neighborhood $U$ of $(0,0)$ where the same strict inequality holds. Take an $\varepsilon > 0$ such that
$$
\{ (\varepsilon , x) : 0 \leq x \leq \varepsilon \} \subset U
$$
Then for all  $ x \, \in  (0, \varepsilon)$ we have
$$
Q(w^+) < 0 < Q(u) = 2h
$$

\item 
To have an upper bound for the normal derivative of $u$ on $\gamma_1$ we need $w^+(0) = 0$ and $w^+(\varepsilon) \geq M$, where 
$$
M = \max_{\{ q \, \in \Omega \, : \, d(q) = \varepsilon\}} \HA(q)
$$
The fact that $w^+(0)=0$ is true by construction of $\psi_{\varepsilon, A}$. To fulfill the second inequality we choose 
$$
A = \dfrac{M}{e^{\varepsilon} - \varepsilon - 1}
$$
We remark that this last choice is possible because $\psi_{\varepsilon, A}'$ is strictly positive  for each $0 \leq x < \varepsilon$, hence $e^{\varepsilon} - \varepsilon - 1$ is strictly greater than $w^+(0) = 0$.
\end{itemize}
We now prove the estimate from below. This will be done using $\HB$ as a punctual lower barrier.\\
Take $z \, \in \, \gamma_1$. From the $r-$admissibility hypothesis made on $\gamma_1$ follows the existence of an horizontal translation $\tau_{z} \left( \HB \right)$ of $\HB$ tangent to $\gamma_1$ in $z$ such that $\tau_z \left( \HB \right) \leq \HA$ and $\tau_z \left( {\HB} \right)_{|\gamma_1} < 0$. Hence by the maximum principle we have
$$
\tau_z \left( \HB \right) \leq u \qquad \mbox{ in } \Omega
$$
\\
From the compactness of $\gamma_1$ together with the regularity of $u$ we get a bound for $||\nabla u||_{C^0(\gamma_1)}$.

\subitem \textit{$C^1$ estimates on $\gamma_2$ }\\
To give an upper barrier we observe that, by the maximum principle, $u$ is below the minimal slice $\Omega \times \{ \HA(\gamma_2) \}$. To give a lower barrier we use a part of a cone, i.e. the graph of an affine function of $\rho$ the hyperbolic distance from $0$. Let's define for $ c >0$ and $k \in \R$ 
\begin{align}
 f(\rho) = c \, \rho + k
\end{align}
One can easily check that, if $\roeucl = \mathrm{tanh} \left( \frac{\rho}{2} \right)$ is the Euclidean distance from $0$, we have
\begin{align}
Q(f) =  \dfrac{c}{\sqrt{ 1 + c^2}} \, \left( \roeucl  + \dfrac{1}{\roeucl \, \lambda(\roeucl)} \right) \label{eq: mean curvature cone}
\end{align}
where $\lambda$ is the conformal factor of the hyperbolic metric in the disc model, see \eqref{eq: lambda}. Note that this curvature is decreasing to the value $\frac{c}{\sqrt{1 + c^2}}$ when $\roeucl$ is increasing to $1$ and this  function of $c$ is monotonically increasing to $1$ when $c$ is increasing to $\infty$.
Now denote by $\rho_2 $ the hyperbolic radius of $\gamma_2$. To define the part of the cone $f$ we are going to use as a barrier we need to choose its domain and its two parameters. The domain will be an annulus of radii $\rho_0 < \rho_2$ and the parameters will be chosen so that on the inner circle the cone is equal to $\HB$ and on the outer it is equal to $\HA$. In other words we look for numbers $\rho_0, c, k$ satisfying the equations
\begin{align}
 c \, \rho_0 + k & = \HB(\rho_0) \\
 c \, \rho_2 + k & = \HA(\rho_2) 
\end{align}

The solutions of these equations, expressed as functions of $\rho_0$, are
\begin{align}
 c & = \dfrac{\HA(\rho_2) - \HB(\rho_0)}{\rho_2  - \rho_0} \label{eq: c} \\
k & = \dfrac{1}{2} \, \Big( \HB(\rho_0)  + \HA(\rho_2) - c \, \left( \rho_0 + \rho_2 \right) \Big)
\end{align}
Hence if we find a $\rho_0$ giving a $c > \frac{2h}{\sqrt{1-4h^2}}$ we are done. Indeed to have a barrier we need a cone with curvature greater than $h$ and taking $c = \frac{2h}{\sqrt{1 - 4h^2}}$ we get a cone with mean curvature monotonically decreasing to $h$. Such a choice for $c$ can be done because, by Theorem \ref{teor: dependance of deriv H in alfa}, the numerator of \eqref{eq: c} goes to a positive number as $\rho_0$ approaches $\rho_2$.

\item \textit{$C^1$ estimates on $\mathrm{int}(\Omega)$}\\
These estimates are a consequence of Theorem 3.1 of Spruck's work \cite{IGSP}. This Theorem gives an interior estimate for $\nabla u$ in terms of the $C^0$ estimate of $u$ and the $C^1$ estimate of $u$ on the boundary, provided $\partial \, \Omega $ is $C^3.$
\end{itemize}
All estimates depending only on translations of $\HA$ and $\HB$ depend only on $h$ and $\gamma_1$ because these two parameters are enough to determine $\HA$ and $\HB$. The only estimate which is not made using these two surfaces is the one made by $\psi_{\varepsilon, A} \circ d$ on $\gamma_1$. It only depends on $\gamma$ because $w^+$ depends on the three parameters $c, \varepsilon$ and $A$, where we recall that
\begin{itemize}
 \item $c$ is the maximum of the geodesic curvature of $\gamma_1$
\item $\varepsilon$ depends on $\gamma_1$
\item $A$ depends on $\gamma_1$ and $\HA$
\end{itemize}

\endproof
\end{teor}
We are now in position to prove the existence Theorem  on compact annuli. We need a slightly modified version of the Method of Continuity allowing a general dependence on the parameter that we use to link the existing solution of a known problem to the one we are looking for. This result can be proved with the standard techniques proposed, for example, in \cite[Chapter 17, Section 2]{GBT}, hence we do not prove it here. The formalism for Fr\'echet differentiability theory is the same as in \cite[Chapter 17, Section 2]{GBT}.
\begin{teor}\label{teor: MdC}$\\$
Let be $0 < \delta < 1$, $\Omega \, \subset \, \R^2$ a compact set and $\mathcal{U} \, \subset \, \CC(\Omega)$ an open set. Consider  a second order  Fr\'echet differentiable operator
$$
\funzione{F}{\mathcal{U}\times [0,1] }{C^{0, \, \delta}(\Omega)}{(u,\sigma)}{F^{\sigma}(u)}
$$
Assume the Dirichlet problem
\begin{align*}
\left \{
\begin{array}{rcccc}
 F(u,1) & = & 0 & \mbox{ in } \quad \Omega & \\
u       & = & \phi  & \mbox{ in } \partial \Omega & \mbox{ for } \quad \phi \, \in C^{2, \, \delta}(\partial \Omega)
\end{array}
\right.
\end{align*}
has solution in $\CC( \Omega)$.\\
If  
\begin{enumerate}[label=\roman*)]
\item $F^\sigma$ is strictly elliptic $\forall \, \sigma \, \in \, [0,1]$
\item $F^{\sigma}_z \leq 0 \quad \forall \, \sigma \, \in [0,1]$
\item $E= \lbrace u \, \in \, \mathcal{U} \, \slash \, \exists \, \sigma \, \in [0,1] \, : \, F(u,\sigma)=0 \quad \mbox{with} \quad u_{|\partial \Omega} = \varphi \rbrace$ is bounded
\item $\bar{E} \subset \mathcal U$
\end{enumerate}
Then the Dirichlet problem
\begin{align*}
\left \{
\begin{array}{cccc}
 F(u,0) & = & 0 & \mbox{ in } \quad \Omega \\
u & = & \phi & \mbox{ in } \quad \partial \Omega
\end{array}
\right.
\end{align*}
has solution in  $\CC( \Omega)$
\end{teor}

We now state our first existence result.
\begin{teor}\label{teor: existence on compact annuli}$\\$
Let be $h \, \in \, (0, \frac{1}{2})$ and $r>0$ and let $\Omega \, \subset \, \hp$ be a compact annulus with boundaries $\gamma_1$ and $\gamma_2$.\\
 Assume $\gamma_1$ is $r-$admissible, horosphere convex and that $\gamma_1$ can be smoothly deformed into a circle, with $r-$admissibility and horosphere convexity preserved along the deformation.\\
 Assume $\gamma_2$ is a circle.\\
Then the Dirichlet problem
 \begin{equation}
  \left \{
\begin{array}{rccr}
  div \left( \dfrac{\nabla u}{\sqrt{1 + {|\nabla u|}^2}} \right) & = & 2 h & \mbox{ in }  \Omega \\
u & = & 0 & \mbox{ in }  \gamma_1 \\ \tag{D} 
u & = & H^h_{\alpha(r)} & \mbox{ in } \gamma_2 
\end{array} \label{eq: D}
\right.
 \end{equation}
has a solution $u \, \in \, C^{2,\delta}(\Omega)$.
\proof
Denote by $D_2$ the disc bounded by $\gamma_2$. We are going to show that we can deform ${\HA}_{ | ( D_2 \setminus  \disco{0}{\rho^h(\alpha)}}$ to a function on $\Omega$ preserving the mean curvature without loosing any regularity. Without loss of generality we can suppose that the deformation assumed in the hypotheses is parametrized by $\sigma \, \in \, [0,1]$ and we write $\gamma^{\sigma}_1$ for the evolution at time $\sigma$ of $\gamma_1$. We also assume $\gamma^0_1 = \gamma_1$ and $\gamma^1_1 = \sfera{0}{\rho^h(\alpha)}$. \\
We define $A^{\sigma}$ to be the annulus whose inner and outer boundaries are respectively $\gamma_1^{\sigma}$ and $\gamma_2$. Then consider $0 < a_1 < a_2$ and $A \, \subset \, \R^2$ the compact annulus with inner boundary the circle of radius $a_1$ and outer boundary the circle of radius $a_2$. Let $\{ \phi^{\sigma} \}_{\sigma \, \in \, [0,1]}$ be a family of smooth orientation preserving embedding of $A$ in $\hp$ such that $\phi^{\sigma}(A) = A^{\sigma}$. The compactness of $[0,1]$ together with the smoothness of the family $\{ \phi^{\sigma} \}_{\sigma}$ yields the existence of a positive constant $C$ which does not depend on $\sigma$ such that:

\begin{align}
 \| \phi_\sigma \|_{\CC(A)} < C
\end{align}
Hence the boundedness of $\| u \|_{\CC(A^\sigma)}$ is equivalent to the boundedness of $\| u \circ \phi^\sigma \|_{\CC(A)}$.\\
We define the operator $F(\sigma, u)$ as follows:
\begin{align*}
 F(u, \sigma) = Q(u \circ \phi^\sigma) - 2h
\end{align*}
For all  $\sigma \, \in [0,1]$ we denote $ u^\sigma = u \circ \phi^\sigma$ and we introduce the family of Dirichlet problems
\begin{equation}\label{probl: dirichlet con sigma}
\left \{
\begin{array}{cccc}
F(u^{\sigma}, \sigma) & = 0  && \mbox{ in } A \nonumber \\ 
u^\sigma & = 0  && \mbox{ in } a_1  \tag{$D^{\sigma}$} \\
u^\sigma & = \HA  && \mbox{ in } a_2 \nonumber 
\end{array}
\right.
\end{equation}
As it is well known, $F$ is Fr\'echet differentiable, for a proof of this fact one can see \cite{BASAMEANCURVATURE}. Moreover $F^{\sigma}$ is constant in $u$ since the mean curvature of a graph does not depend on the height where it is evaluated, hence we have $F_z \leq 0$.\\
We now state bounds for $|| u ||_{C^1(A_{\sigma})}$ independent of $\sigma$. All the estimates given in terms of $\HA$ and $\HB$ do not depend on $\sigma$ because $r-$admissibility is preserved in $\sigma$ and hence neither $\HA$ nor $\HB$ depend on the parameter. Redefining $c$ in \eqref{eq: max geod curvature along gamma_1} by
$$
c = \max_{\begin{array}{c} \sigma \, \in \, [0,1]\\ y_0 \, \in \, \gamma_1^\sigma \end{array}} \, k_{g}(\gamma_1^\sigma)(y_0)
$$
we obtain an estimate of $\nabla u^\sigma$ on ${\gamma_1}^\sigma$ independent from $\sigma$.
\endproof
\end{teor}

For examples of curves satisfying to the hypotheses of the existence Theorem we refer to the remarks following the Definition of $r-$admissibility.\\
We end the section with the main geometric result of this work: the existence of non rotational vertical ends. The existence of our cmc non rotational ends follows from our existence result on the compact case. Indeed we consider a sequence of compact annuli diverging to the exterior domain on which we are building the end. Then, by means of Theorem \ref{teor: existence on compact annuli}, we obtain a sequence of compact cmc graphs. Each of these graphs on the outer boundary coincides with $\HA$, hence the limit of the graphs will have the same asymptotic behavior of $\HA$. 
\begin{teor}\label{teor: existence on exterior}$\\$
 Let $h \, \in \, (0, \frac{1}{2})$ and $r>0$. Let $\Omega$ be the complement of a compact domain of $\hp$ with boundary a Jordan smooth curve $\gamma_1$. Assume $\gamma_1$ is $r-$admissible, horosphere convex and that it can be smoothly deformed into a circle, with $r-$admissibility and horosphere convexity preserved along the deformation.\\ Then the following Dirichlet problem

\begin{align}
 \left \{
\begin{array}{rccr}
 \dfrac{1}{2} \, div \left( \dfrac{\nabla u}{\sqrt{1 + {|\nabla u|}^2}} \right) & = & h & \mbox{ in }  \Omega \\
u & = & 0 & \mbox{ in }  \gamma_1 \\ \tag{E}
\lim_{|z|_{\hp} \to + \infty} u(z) & = & +\infty & 
\end{array}
\right.
\end{align}
has a solution $u \, \in \, \CC(\Omega)$ for some $0 < \delta < 1$.

\proof
To prove the Theorem we proceed as follows. We consider a sequence of compact annuli $\Omega_n$ converging to $\Omega$, we build a sequence of $\CC(\Omega)$ solutions by mean of Theorem \ref{teor: existence on compact annuli} and then we prove convergence in $\CC(\Omega)$.\\
To accomplish the first step we consider $ \{ \rho_n \}_{n \, \in \N}$ a sequence of positive reals monotonically diverging to $+ \infty$. We define $\gamma_n$ to be the circle $\sfera{0}{\rho_n}$. Then we define $\Omega_n$ to be the annulus whose inner and outer boundary are respectively $\gamma_1$ and $\gamma_n$. We introduce the Dirichlet problem
\begin{align}
 \left \{
\begin{array}{cccc}
 div \left( \dfrac{\nabla u_n}{\sqrt{1 + {|\nabla u|}^2}} \right) & = 2\, h  && \mbox{ in } \Omega_n \nonumber \\ 
u_n & = 0  && \mbox{ in } \gamma_1  \\
u_n & = \HA  && \mbox{ in } \gamma_n  \tag{$D_n$} 
\end{array}\label{eq: Dirchlet ennesimo}
\right.
\end{align}
For all $n \, \in \, \N$, Theorem \ref{teor: existence on compact annuli} implies existence $u_n \, \in \, \CC(\Omega_n)$ a solution of problem \eqref{eq: Dirchlet ennesimo}. Moreover Theorem \ref{teor: estimates on annuli} and following remarks give an estimate for $\| u_n \|_{\CC(\Omega)}$ which depends on $n$ only because of the $C^1$ estimate on $\gamma_n$ . 
 Hence we only have to prove the independence on $n$ of the estimates given by Theorem \ref{teor: estimates on annuli}. \\
Now let's give a lower barrier for $u_n$ on $\gamma_n$ so that the $C^1$ estimate is not dependent on $n$. We are going to exhibit a sequence barriers of the same kind we used in the proof oh Theorem \ref{teor: existence on compact annuli}, which means we are going to introduce a sequence of cones $\{ f \}_n$ whose normal derivatives are bounded in $n$ and each $f_n$ is a lower barrier for the gradient of $u_n$ on $\gamma_n$.\\
Denote by $\rho_n$ the radius of $\gamma_n$. We recall that, by hypothesis of $r-$admissibility on $\Omega$, we associate to each $\Omega_n$ the same two surfaces $\HA$ and $\HB$. This is because the choice of $\HA$ and $\HB$ is made on $\gamma_1$ (see the Definition of $r-$admissibility \ref{def: curva r-ammissibile}).\\
For $\varepsilon \, \in \, (0, \rho_n - r)$ we denote
\begin{align*}
 t_{ \varepsilon} & = \HB(\rho_n - \varepsilon) \\
t_n & = \HA(\rho_n)
\end{align*}
If $n$ is large exists $\varepsilon \, \in \, (0, \rho_n - r)$ such that the function
\begin{align*}
 \funzione{f_{n, \varepsilon}}{A(\rho_n - \varepsilon, \rho_n)}{\R}{\rho}{c_{n, \varepsilon} \, \rho + k_{n,\varepsilon}}
\end{align*}
with
\begin{align*}
 c_{n, \varepsilon} & = \dfrac{t_n - t_{\varepsilon}}{\varepsilon}\\
k_{n,\varepsilon} & = t_n - c_{n, \varepsilon} \, \rho_n
\end{align*}
is a lower barrier for $u_n$ on $\gamma_n$. Moreover for each fixed $\varepsilon > 0$ the sequence $\{ c_{n, \varepsilon}\}_n$ is bounded, hence each $f_{n,\varepsilon}$ has normal derivative bounded in $n$.\\
The boundary data are assumed, i.e. $f_n(\rho_n - \varepsilon ) = t_{\varepsilon}$ and $f_n(\rho_n) = t_n$. Now we check that exists $\varepsilon$ such that for all $n$ we have
\begin{align*}
 Q(f_n) > 2h
\end{align*}
Using the asymptotic expansion \eqref{eq: sviluppo asintotico con alpha}  of $\HA$ and $\HB$ we get 
\begin{align}
c_{n, \varepsilon} = \dfrac{k_{\alpha} - k_{\beta}}{\varepsilon} + \dfrac{2h}{\sqrt{1 - 4h^2}} \, \varepsilon + o(e^{- \rho_n}) \label{eq: espressione delle c_n} 
\end{align}
where the first term is positive because of Theorem \ref{teor: dependance of deriv H in alfa}. Recalling the value of the mean curvature of a cone given in \eqref{eq: mean curvature cone} one can see that $Q(f_n)$ decreases to $\frac{c_{n, \varepsilon}}{\sqrt{1 + c_{n, \varepsilon}^2}}$. Hence if
\begin{align}
 \dfrac{c_{n,\varepsilon}}{\sqrt{1 + c_{n, \varepsilon}^2 }} \geq 2h \label{eq: curvatura cono grande abbastanza}
\end{align}
then for all $n \, \in \, \N$
$$
Q(f_n) > 2h
$$
Being that the mean curvature of $f_{n, \varepsilon}$  increases to $1$ as $c_{n,\varepsilon}$ increases to $+\infty$ and being  $2 h < 1$, it is clear that to satisfy \eqref{eq: curvatura cono grande abbastanza} it is enough to choose $\varepsilon >0$ such that
$$
\varepsilon \leq \dfrac{k^h_{\alpha} - k^h_{\beta}}{\frac{2h}{\sqrt{1- 4h^2}}}  
$$ 
We have proved that $f_n$ is a lower barrier for $u_n$ and hence it gives a lower bound for the normal derivative of $u_n$ on $\gamma_n$. To prove that the normal derivative of each $f_n$ on $\gamma_n$ is bounded in $n$ it is enough to recall that this derivative coincides with $c_{n, \varepsilon}$ and, by equation \eqref{eq: espressione delle c_n}, the sequence $\{ c_{n, \varepsilon} \}_n$ converges and hence is bounded.

We have proved that $\{ u_n \}_n$ is a sequence of $C^{2,\delta}$ functions which is bounded and equicontinuous on the compacts of $\Omega$. Hence we can apply Ascoli - Arzel\`a Theorem and, up to passing to a subsequence, we have the convergence in $\CC(\Omega)$.
\endproof
\end{teor}

\newpage
\bibliography{biblio}
\bibliographystyle{amsplain}

\end{document}